\documentclass[11pt]{amsart}

\usepackage{amssymb}
\usepackage{amsmath}
\usepackage{amsthm}
\usepackage{tikz}
\usetikzlibrary{arrows}
\usepackage[all]{xy}
\xyoption{2cell}
\usepackage{graphicx}
\usepackage{hyperref}

\DeclareMathOperator{\Trop}{Trop}
\DeclareMathOperator{\orb}{orb}
\DeclareMathOperator{\Div}{Div}

\DeclareMathOperator{\Pic}{Pic}

\DeclareMathOperator{\Spec}{Spec}
\DeclareMathOperator{\val}{val}

\DeclareMathOperator{\Nef}{Nef}

\DeclareMathOperator{\Hom}{Hom}

\theoremstyle{plain}
\newtheorem{theorem}{Theorem}[section]
\newtheorem{lemma}[theorem]{Lemma}
\newtheorem{proposition}[theorem]{Proposition}
\newtheorem{corollary}[theorem]{Corollary}

\theoremstyle{definition}
\newtheorem{definition}[theorem]{Definition}

\newtheorem{construction}[theorem]{Construction}

\theoremstyle{remark}
\newtheorem{remark}[theorem]{Remark}

\begin{document}

\title{Numerical tropical line bundles and toric b-divisors}

\author{Carla Novelli  and Stefano Urbinati}
\date{}
\address{Carla Novelli - Universit\`a degli Studi di Torino}
\email{carla.novelli@unito.it} 
\address{Stefano Urbinati - Universit\`a degli Studi di Udine} \email{stefano.urbinati@uniud.it}
\maketitle

\begin{abstract}
We study the relationship between line bundles on tropical compactifications of a very affine variety $Y$ and toric b-divisors on the associated tropical variety $\Trop(Y)$. By focusing on numerical equivalence classes, we construct a natural injective map from the group of numerical tropical line bundles on $Y$ to the space of toric b-divisors modulo linear equivalence. Moreover, we show that this map restricts to a bijection between the tropical nef cone of $Y$ and the set of toric b-divisors that are b-Cartier and tropically nef. This provides a higher-dimensional generalization of Baker's specialization for curves and clarifies the birational nature of tropical line bundles. We also discuss the kernel of the map from line bundles to numerical tropical line bundles, which encodes the continuous moduli lost in tropicalization.
\end{abstract}

\section{Setup and assumptions}
\label{sec:setup}

Throughout this paper, we work with a very affine variety $Y$, i.e., a closed subvariety of an algebraic torus $T_N \cong (\mathbb{C}^*)^n$. We assume that $Y$ is \emph{sch\"on} in the sense of Tevelev \cite{Tev-CptSubvar}, meaning that it admits tropical compactifications with proper intersections. This assumption is essential for the following reasons:

\begin{enumerate}
\item It guarantees the existence of a tropical fan $\Sigma$ with $|\Sigma| = \Trop(Y)$ such that the closure $\overline{Y}_\Sigma$ in $\mathbb{P}_\Sigma$ is a tropical compactification. After refining $\Sigma$ if necessary, $\overline{Y}_\Sigma$ can be made smooth and its boundary a reduced simple normal crossing divisor.
\item For any tropical fan $\Sigma$, the intersection of $\overline{Y}_\Sigma$ with any torus orbit $\orb(\sigma)$ is either empty or of the expected dimension $\dim Y - \dim \sigma$.
\item The numerical intersection theory on $\overline{Y}_\Sigma$ is controlled by torus-invariant cycles.
\end{enumerate}

All constructions and results in this paper depend on these properties.

\section{Introduction}
\label{intro}

The study of line bundles on tropical varieties has been a central theme in tropical geometry, with connections to arithmetic, algebraic geometry, and combinatorics. A fundamental problem is to understand how algebraic line bundles on a variety $Y$ induce tropical line bundles on its tropicalization $\Trop(Y)$. 

In this paper, we revisit this problem from the perspective of birational geometry and valuation theory. Our starting point is the observation that tropical compactifications, introduced by Tevelev \cite{Tev-CptSubvar}, provide a natural framework for extending line bundles from a very affine variety $Y$ to its compactifications in toric varieties. The non-uniqueness of such extensions reflects the birational nature of the problem.

We interpret this non-uniqueness using the language of b-divisors (short for \emph{birational divisors}), i.e., collections of divisors on all birational models, compatible under pullback. This approach, developed by Shokurov and others, allows us to work consistently across different tropical compactifications. Specifically, we show that to every line bundle on a tropical compactification of $Y$, one can associate a toric b-divisor. This construction is compatible with refinements of the tropical fan, yielding a well-defined map from tropical line bundles to toric b-divisors.

However, as recognized in the literature on tropical curves \cite{Baker}, this map cannot be injective when considering isomorphism classes of line bundles, since tropicalization forgets continuous parameters. We therefore focus on \emph{numerical equivalence classes}, which capture precisely the discrete data preserved by tropicalization. Our main result (Theorem \ref{thm:main}) establishes that the induced map from the group of numerical tropical line bundles to the space of toric b-divisors modulo linear equivalence is injective. Furthermore, we characterize its image when restricted to the nef cone: a numerical tropical line bundle is nef if and only if its associated b-divisor is b-Cartier and tropically nef, giving a bijection between the tropical nef cone and the set of such b-divisors. This provides a higher-dimensional generalization of Baker's specialization map for curves.

The paper is organized as follows. Section 2 reviews background on toric and tropical geometry, tropical compactifications, and Minkowski weights. Section 3 introduces the necessary notions from valuation theory and b-divisors, focusing on the aspects relevant to the tropical setting. Section 4 contains the main construction and results. Section 5 provides illustrative examples, and Section 6 discusses connections with existing work and open questions.

\section{Preliminaries on toric and tropical geometry}

\subsection{Toric varieties and divisors}
We briefly recall standard notions from toric geometry \cite{toric}. Let $N \cong \mathbb{Z}^n$ be a lattice, $M = \Hom(N,\mathbb{Z})$ its dual, and $\Sigma$ a fan in $N_{\mathbb{R}} = N \otimes_{\mathbb{Z}} \mathbb{R}$. The associated toric variety is denoted $\mathbb{P}_{\Sigma}$, with dense torus $T_N \cong (\mathbb{C}^*)^n$.

Each ray $\rho \in \Sigma(1)$ corresponds to a torus-invariant prime divisor $D_{\rho}$. A $T_N$-invariant Weil divisor is of the form $D = \sum_{\rho \in \Sigma(1)} a_{\rho} D_{\rho}$. Such a divisor is Cartier if for each maximal cone $\sigma \in \Sigma(n)$ there exists $m_{\sigma} \in M$ with $D|_{U_{\sigma}} = \operatorname{div}(\chi^{m_{\sigma}})$. The polyhedron of $D$ is
\[
P_D = \{ m \in M_{\mathbb{R}} \mid \langle m, u_{\rho} \rangle \geq -a_{\rho} \text{ for all } \rho \in \Sigma(1) \},
\]
where $u_{\rho}$ is the primitive generator of $\rho$.

\subsection{Very affine varieties and tropicalization}
A \emph{very affine variety} is a closed subvariety of an algebraic torus. We will assume that $Y \subseteq T$ is very affine and irreducible, so that the group $\mathcal{O}^*(Y)/\mathbb{C}^*$ is a lattice $M_Y$, and the embedding $Y \hookrightarrow T_{M_Y} := \Spec \mathbb{C}[M_Y]$ is universal up to the action of $\mathbb{C}^*$ on the torus; different choices give embeddings that differ by a torus automorphism. We call $T_{M_Y}$ the \emph{intrinsic torus} of $Y$.

For a very affine variety $Y \subseteq (\mathbb{C}^*)^n$, the \emph{tropicalization} $\Trop(Y)$ is defined as follows. Extend the ground field to the field of Puiseux series $\mathbb{C}\{\!\{t\}\!\}$ with the usual valuation $\val$; then $Y$ gives rise to a variety $Y^{\mathrm{an}}$ over this field. The tropicalization $\Trop(Y)$ is the closure in $\mathbb{R}^n$ of the set
\[
\{ (\val(x_1), \dots, \val(x_n)) \mid (x_1,\dots,x_n) \in Y^{\mathrm{an}}(\mathbb{C}\{\!\{t\}\!\}) \}.
\]
It is a rational polyhedral set of pure dimension $\dim Y$. For a sch\"on variety, $\Trop(Y)$ is the support of a rational polyhedral fan. Moreover, there exists a fan $\Sigma$ (called a \emph{tropical fan}) with $|\Sigma| = \Trop(Y)$ such that $(\overline{Y}_\Sigma, \mathbb{P}_\Sigma)$ is a tropical compactification; every refinement of $\Sigma$ is again a tropical fan. In this paper we will work with such a tropical fan $\Sigma$ (or a refinement thereof), and we will refer to $|\Sigma|$ as $\Trop(Y)$.

\subsection{Schön varieties}
\label{subsec:schon-def}

A key concept in the theory of tropical compactifications is the notion of a \emph{schön} variety. Introduced by Tevelev~\cite{Tev-CptSubvar}, this property guarantees that the geometry of the compactification is particularly well-behaved.

\begin{definition}\label{def:schon}
Let $Y \subseteq T$ be a very affine variety. We say that $Y$ is \emph{schön} if there exists a unimodular fan $\Sigma$ with support $|\Sigma| = \Trop(Y)$ such that the open strata $Y^\sigma := \overline{Y}_\Sigma \cap \orb(\sigma)$ are smooth for every cone $\sigma \in \Sigma$ \cite[Definition 1.1]{Aksnes2025}.
\end{definition}

A fundamental result, which follows from the work of Tevelev and Hacking, is that this property is independent of the chosen fan.

\begin{theorem}[{\cite[Theorem 4.4]{Aksnes2025}, \cite{LuxtonQu2011}}] \label{thm:schon}
If a very affine variety $Y$ is schön, then there exists a tropical fan $\Sigma$ (i.e., a fan such that $(\overline{Y}_\Sigma, \mathbb{P}_\Sigma)$ is a tropical compactification) with $|\Sigma| = \Trop(Y)$. Moreover, for every such tropical fan $\Sigma$ (and in particular for every refinement of $\Sigma$ that is still a tropical fan), the following properties hold:
\begin{enumerate}
    \item The open strata $Y^\sigma = \overline{Y}_\Sigma \cap \orb(\sigma)$ are smooth for all $\sigma \in \Sigma$.
    \item After a possible refinement of $\Sigma$, the compactification $\overline{Y}_\Sigma$ can be made smooth and the boundary $\overline{Y}_\Sigma \setminus Y$ is a reduced simple normal crossing divisor.
    \item The intersection of $\overline{Y}_\Sigma$ with any torus orbit $\orb(\sigma)$ is either empty or of the expected dimension $\dim Y - \dim \sigma$.
\end{enumerate}
\end{theorem}

\begin{remark}
The schön condition is not overly restrictive; in characteristic zero, every algebraic variety contains a Zariski-open very affine subvariety that is schön~\cite{LuxtonQu2011}. This makes the class of schön varieties a rich and important testing ground for non-archimedean and tropical methods.
\end{remark}

\subsection{Tropical compactifications}
\begin{definition}[\cite{Tev-CptSubvar}]\label{def:tropical-compactification}
Let $Y \subset T$ be a very affine variety. Let $\Sigma$ be a fan in $N_{\mathbb{R}}$ such that $|\Sigma| = \Trop(Y)$. Let $\overline{Y}$ be the closure of $Y$ in the toric variety $\mathbb{P}_{\Sigma}$. The pair $(\overline{Y}, \mathbb{P}_{\Sigma})$ is a \emph{tropical compactification} if the multiplication map
\[
\overline{Y} \times T \longrightarrow \mathbb{P}_{\Sigma}
\]
is flat and surjective. In this case, $\Sigma$ is called a \emph{tropical fan} for $Y$.
\end{definition}

Tropical compactifications have favorable properties: the boundary $\overline{Y} \setminus Y$ is a reduced divisor with simple normal crossings (after suitable refinement), and $\overline{Y}$ intersects each torus orbit properly:
\begin{equation}\label{eq:proper_int}
\dim(\orb(\sigma) \cap \overline{Y}) = \dim Y - \dim \sigma \quad \text{for all } \sigma \in \Sigma.
\end{equation}

By \cite[Proposition 2.3]{Tev-CptSubvar}, any refinement of a tropical fan that remains a fan is again tropical. Thus, tropical fans form a directed system under refinement. In particular, by refining the fan we can resolve singularities of the toric variety, making it smooth.

\subsection{Construction of the toric fan from tropicalization}
\label{subsec:fan-construction}
We recall a construction from \cite{KS} that builds a toric scheme over a DVR starting from a polyhedral subdivision of $\Trop(Y)$. Although $Y$ is defined over $\mathbb{C}$, one can consider the formal power series ring $\mathbb{C}[[t]]$ and its fraction field $\mathbb{C}((t))$ to obtain a one‑parameter degeneration. More precisely, let $\Sigma$ be a rational polyhedral fan with support $|\Sigma| = \Trop(Y)$. Consider $\Sigma \times \{1\} \subseteq N_{\mathbb{R}} \times \mathbb{R}_{\ge 0}$ and let $\tilde{\Sigma}$ be the fan of cones over the facets of $\Sigma \times \{1\}$. The toric variety $\mathbb{P}_{\tilde{\Sigma}}$ comes with a projection $\mathbb{P}_{\tilde{\Sigma}} \to \mathbb{A}^1$ induced by the second factor. For the DVR $\mathcal{O} = \mathbb{C}[[t]]$ with uniformizer $t$, the base change $\mathbb{P}_{\Sigma} := \mathbb{P}_{\tilde{\Sigma}} \times_{\Spec \mathbb{Z}[u]} \Spec \mathcal{O}$ gives a toric scheme whose generic fiber is the toric variety associated to $\Sigma$ (now considered over $\mathbb{C}((t))$). This construction shows how the tropical fan naturally arises from a formal degeneration; the original variety $Y$ remains defined over $\mathbb{C}$ and is recovered in the special fiber after suitable modifications.

\begin{remark}
If one works with the trivial valuation, $\mathbb{P}_{\Sigma}$ can be compactified to a complete toric variety by adding cones corresponding to the missing directions; the resulting compactification has a boundary of codimension at least $2$.
\end{remark}

\subsection{Minkowski weights}
A tropical compactification induces a \emph{Minkowski weight} on the fan $\Sigma$.
\begin{definition}
\label{def:minkowski-weight}
Let $\Sigma$ be a pure $d$-dimensional fan. A function $c \colon \Sigma(d) \to \mathbb{Z}$ is a \emph{Minkowski weight of dimension $d$} if for every $\tau \in \Sigma(d-1)$,
\[
\sum_{\substack{\sigma \in \Sigma(d) \\ \sigma \supset \tau}} c(\sigma) \, v_{\sigma/\tau} = 0 \quad \text{in } N/N_{\tau},
\]
where $v_{\sigma/\tau}$ is a primitive generator of $\sigma$ modulo $\tau$.
\end{definition}

Given a tropical compactification $\overline{Y} \subset \mathbb{P}_{\Sigma}$ of dimension $k$, the associated Minkowski weight $c_Y$ is defined by $c_Y(\sigma) = \deg([\overline{Y}] \cdot [V(\sigma)])$ for $\sigma \in \Sigma(k)$. This weight is balanced and represents the tropicalization of $Y$.

\begin{lemma}[{\cite[Lemma 6.2]{Ka-TropIntToric}}]\label{lem:strata-eq}
Let $\mathcal{L}_1, \mathcal{L}_2 \in \Pic(\mathbb{P}_{\Sigma})$ be two line bundles on the toric variety $\mathbb{P}_{\Sigma}$. Assume that for every cone $\tau \in \Sigma(k-1)$ such that the curve $C_\tau := \overline{Y}_{\Sigma} \cap V(\tau)$ is non‑empty, we have
\[
\deg\bigl(\mathcal{L}_1|_{C_\tau}\bigr) = \deg\bigl(\mathcal{L}_2|_{C_\tau}\bigr).
\]
Then the restrictions $\mathcal{L}_1|_{\overline{Y}_{\Sigma}}$ and $\mathcal{L}_2|_{\overline{Y}_{\Sigma}}$ are numerically equivalent on $\overline{Y}_{\Sigma}$.
\end{lemma}


This lemma will be used to prove the injectivity of the main map.

\section{Valuations and b-divisors in the toric setting}

We now recall the necessary notions from valuation theory and b-divisors, following \cite{RU}, with a focus on the toric setting.

\subsection{Toric valuations and the Berkovich analytification}
Let $T \cong (\mathbb{C}^*)^n$ be an algebraic torus with character lattice $M$. The Berkovich analytification $T^{\mathrm{an}}$ of $T$ (with respect to the trivial valuation on $\mathbb{C}$) contains all multiplicative seminorms on $\mathbb{C}[M]$ extending the trivial norm on $\mathbb{C}$.

A point $\xi \in T^{\mathrm{an}}$ can be described as a pair $(V, v)$ where $V \subseteq T$ is an irreducible subvariety and $v$ is a real valuation on $\mathbb{C}(V)$ centered on $T$. Of particular interest are \emph{toric valuations}, which are quasi-monomial valuations that become monomial after a toric birational modification.

\begin{definition}
A valuation $v$ on $\mathbb{C}(T)$ is called \emph{toric} if there exists a toric variety $\mathbb{P}_{\Sigma}$ (with $\Sigma$ a fan in $N_{\mathbb{R}}$) and a point $p$ in the torus orbit corresponding to a cone $\sigma \in \Sigma$ such that $v$ is monomial with respect to local coordinates at $p$ that are monomials in the toric coordinates.
\end{definition}

Toric valuations are precisely the quasi-monomial valuations for which the associated graded algebra is a monomial algebra. They correspond to points in the tropicalization of $T$, which is $\mathbb{R}^n$. Moreover, for a schön variety $Y$, the restriction of such valuations to $\mathbb{C}(Y)$ yields points in $\Trop(Y)$ interpreted as a fan.

\subsection{b-divisors on toric Zariski--Riemann spaces}

Let $\Sigma_0$ be a tropical fan for $Y$ as in Definition~\ref{def:num-trop}. The \emph{toric Zariski--Riemann space} $\mathfrak{X}_{\mathrm{tor}}$ is defined as the inverse limit of all toric varieties $\mathbb{P}_{\Sigma}$ where $\Sigma$ is a refinement of $\Sigma_0$:
\[
\mathfrak{X}_{\mathrm{tor}} = \varprojlim_{\Sigma \geq \Sigma_0} \mathbb{P}_{\Sigma}.
\]
A \emph{toric b-divisor} is an element of $\Div(\mathfrak{X}_{\mathrm{tor}}) = \varprojlim_{\Sigma} \Div(\mathbb{P}_{\Sigma})$.

Given a toric valuation $v$ (or equivalently a point $\xi \in T^{\mathrm{an}}$ corresponding to a toric valuation), we can associate a toric b-divisor $D_{\xi}$ following the construction in \cite{RU}. For each toric model $\mathbb{P}_{\Sigma'}$, consider the valuation ideal sheaf $\mathfrak{a}_{\xi, m}^{\Sigma'} = \{ f \in \mathcal{O}_{\mathbb{P}_{\Sigma'}} \mid v(f) \geq m \}$. Let $Z(\mathfrak{a}_{\xi, m}^{\Sigma'})$ be the Cartier b-divisor defined by the normalized blowup of $\mathfrak{a}_{\xi, m}^{\Sigma'}$. Then
\[
D_{\xi} = \lim_{m \to \infty} \frac{1}{m} Z(\mathfrak{a}_{\xi, m}^{\Sigma_0}),
\]
where $\Sigma_0$ is some initial fan. The limit exists as a nef b-divisor and is independent of the choice of $\Sigma_0$ up to refinement.

\begin{proposition}[{\cite[Proposition 2.13]{RU}}]
For a toric valuation $v$, the b-divisor $D_{\xi}$ is toric and can be described explicitly in terms of the weight vector defining $v$. Moreover, $D_{\xi}$ is nef and anti-effective.
\end{proposition}

\subsection{Tropical nefness and b-Cartier property}
In the toric context, we need adapted notions of nefness and the Cartier property for b-divisors.

\begin{definition}\label{def:tropical-nef}
Let $\mathbf{D} = (D_{\Sigma'})$ be a toric b-divisor on $\mathfrak{X}_{\mathrm{tor}}$.
We say that $\mathbf{D}$ is \emph{tropically nef} if there exists a refinement $\Sigma'$ (of the fan $\Sigma_0$ with respect to which $\mathfrak{X}_{\mathrm{tor}}$ is defined) such that the piecewise linear function $\phi_{\Sigma'} : |\Sigma'| \to \mathbb{R}$ associated to $D_{\Sigma'}$ (uniquely determined up to a global linear function) has non‑negative slopes across every codimension‑one wall. By the compatibility condition, this property is independent of the chosen refinement.
\end{definition}

\begin{definition}
\label{def:b-cartier}
A toric b-divisor $\mathbf{D}$ is \emph{b-Cartier} if there exists a refinement $\Sigma'$ such that for every further refinement $\Sigma'' \geq \Sigma'$, we have $D_{\Sigma''} = \pi_{\Sigma''/\Sigma'}^* D_{\Sigma'}$, where $\pi_{\Sigma''/\Sigma'} \colon \mathbb{P}_{\Sigma''} \to \mathbb{P}_{\Sigma'}$ is the toric morphism.
\end{definition}

These properties will characterize the image of the map from tropical line bundles to toric b-divisors.

\section{Numerical tropical line bundles}

\begin{definition}\label{def:num-trop}
Let $Y$ be a sch\"on very affine variety. Fix a tropical fan $\Sigma_0$ for $Y$ (which exists by Theorem~\ref{thm:schon}). Consider the directed set of all refinements $\Sigma$ of $\Sigma_0$ such that $\Sigma$ is again a tropical fan for $Y$. For each such $\Sigma$, let $N^1(\overline{Y}_{\Sigma})$ be the N\'eron--Severi group of the tropical compactification $\overline{Y}_{\Sigma}$. If $\Sigma'$ is a refinement of $\Sigma$, the natural birational morphism $\overline{Y}_{\Sigma'} \to \overline{Y}_{\Sigma}$ induces a pullback map $N^1(\overline{Y}_{\Sigma}) \to N^1(\overline{Y}_{\Sigma'})$. The group of \emph{numerical tropical line bundles} on $Y$ is defined as the direct limit
\[
N^1_{\mathrm{trop}}(Y) := \varinjlim_{\Sigma} N^1(\overline{Y}_{\Sigma}).
\]
This definition is independent of the choice of the initial fan $\Sigma_0$ up to canonical isomorphism.
\end{definition}

\begin{remark}
The definition of $N^1_{\mathrm{trop}}(Y)$ is independent of the choice of the initial tropical fan $\Sigma_0$.
Indeed, given two tropical fans $\Sigma_0$ and $\Sigma_0'$ supported on $\Trop(Y)$, there exists a common refinement
$\Sigma''$ refining both. It follows that the corresponding directed systems are cofinal, and hence their direct limits
are canonically isomorphic.
\end{remark}

\begin{remark}
An element of $N^1_{\mathrm{trop}}(Y)$ can be represented by a pair $(\Sigma, [\mathcal{L}])$ where $[\mathcal{L}] \in N^1(\overline{Y}_{\Sigma})$, with $(\Sigma, [\mathcal{L}]) \sim (\Sigma', [\mathcal{L}'])$ if their pullbacks to a common refinement coincide numerically.
\end{remark}

\subsection{Tropical weights of line bundles}
\begin{construction}\label{const:weight}
Let $\mathcal{L} \in \Pic(\overline{Y}_{\Sigma})$. For each cone $\tau \in \Sigma(k-1)$ (where $k = \dim Y$), let $V(\tau)$ be the corresponding torus-invariant subvariety of $\mathbb{P}_{\Sigma}$. By \eqref{eq:proper_int}, $C_{\tau} := \overline{Y}_{\Sigma} \cap V(\tau)$ is a curve (possibly reducible). Define the \emph{tropical weight} of $\mathcal{L}$ along $\tau$ as
\[
w_{\mathcal{L}}(\tau) := \deg\left( \mathcal{L}|_{C_{\tau}} \right).
\]
If $\dim \tau' < k-1$, then $C_{\tau''}$ is contracted by $\pi_{\Sigma''/\Sigma'}$, hence
\[
\deg(\mathcal{L}''|_{C_{\tau''}}) = 0.
\]
This is compatible with the fact that the corresponding slope of the piecewise linear function is zero,
since pullback preserves linearity on cones of lower dimension.
\end{construction}

\begin{lemma}\label{lem:balancing}
The function $w_{\mathcal{L}} \colon \Sigma(k-1) \to \mathbb{Z}$ defined in Construction \ref{const:weight}
satisfies the balancing condition: for every $\gamma \in \Sigma(k-2)$,
\[
\sum_{\substack{\tau \in \Sigma(k-1) \\ \tau \supset \gamma}} w_{\mathcal{L}}(\tau) \cdot v_{\tau/\gamma} = 0 \quad \text{in } N/N_{\gamma},
\]
where $v_{\tau/\gamma}$ is a primitive generator of $\tau$ modulo $\gamma$.
Hence $w_{\mathcal{L}}$ is a Minkowski weight of codimension one on $\Trop(Y)$.
\end{lemma}

\begin{proof}
Let $\mathcal{L}$ be a line bundle on $\overline{Y}_{\Sigma}$. For each $\tau \in \Sigma(k-1)$, recall that
\[
w_{\mathcal{L}}(\tau) = \deg(\mathcal{L}|_{C_\tau}),
\quad \text{where } C_\tau = \overline{Y}_{\Sigma} \cap V(\tau).
\]

We reinterpret these numbers intersection-theoretically. Since $Y$ is schön and $\Sigma$ is a tropical fan, all toric strata $V(\sigma)$ intersect $\overline{Y}_{\Sigma}$ properly. In particular, for each $\tau \in \Sigma(k-1)$,
\[
w_{\mathcal{L}}(\tau) = \deg\bigl(c_1(\mathcal{L}) \cdot [\overline{Y}_{\Sigma}] \cdot [V(\tau)]\bigr).
\]

Fix $\gamma \in \Sigma(k-2)$ and consider the surface
\[
S_\gamma := \overline{Y}_{\Sigma} \cap V(\gamma).
\]
By the schön property, this is a proper surface in the toric variety $V(\gamma)$, and its boundary inside $V(\gamma)$ is given by
\[
\partial S_\gamma = \bigcup_{\tau \supset \gamma} C_\tau.
\]

Let $u \in M := \operatorname{Hom}(N,\mathbb{Z})$ be any character that vanishes on $N_\gamma$, and consider the corresponding monomial function $\chi^u$ on the torus. Its divisor on $V(\gamma)$ is
\[
\operatorname{div}(\chi^u)|_{V(\gamma)} = \sum_{\tau \supset \gamma} \langle u, v_{\tau/\gamma} \rangle \, V(\tau),
\]
where $v_{\tau/\gamma}$ denotes the primitive generator of $\tau$ modulo $\gamma$.

Intersecting with $c_1(\mathcal{L}) \cdot [\overline{Y}_{\Sigma}]$ and using the projection formula, we obtain
\[
0
= \deg\bigl(c_1(\mathcal{L}) \cdot [\overline{Y}_{\Sigma}] \cdot \operatorname{div}(\chi^u)\bigr)
= \sum_{\tau \supset \gamma} \langle u, v_{\tau/\gamma} \rangle \, w_{\mathcal{L}}(\tau).
\]

Since this holds for all $u \in M$ vanishing on $N_\gamma$, it follows that
\[
\sum_{\tau \supset \gamma} w_{\mathcal{L}}(\tau) \, v_{\tau/\gamma} = 0
\quad \text{in } N/N_\gamma.
\]

This proves the balancing condition, hence $w_{\mathcal{L}}$ is a Minkowski weight of codimension one.
\end{proof}

\subsection{Associated toric b-divisor}
\begin{construction}\label{const:bdiv}
Given $[\mathcal{L}] \in N^1(\overline{Y}_{\Sigma})$ represented by a line bundle $\mathcal{L}$, we define a toric b-divisor $\mathbf{D}_{[\mathcal{L}]}$ on $\mathfrak{X}_{\mathrm{tor}}$ as follows.

For each refinement $\Sigma'$ of $\Sigma$, let $\pi_{\Sigma'} \colon \mathbb{P}_{\Sigma'} \to \mathbb{P}_{\Sigma}$ be the corresponding toric morphism, and let $\overline{Y}_{\Sigma'}$ be the closure of $Y$ in $\mathbb{P}_{\Sigma'}$ (the strict transform of $\overline{Y}_{\Sigma}$). Let $\mathcal{L}' = \pi_{\Sigma'}^* \mathcal{L}|_{\overline{Y}_{\Sigma'}}$.

The weights $w_{\mathcal{L}'}(\tau')$ (defined as in Construction \ref{const:weight}) form a balanced collection on $\Sigma'$ by Lemma \ref{lem:balancing}. By the standard correspondence between Minkowski weights of codimension one and piecewise linear functions on a fan \cite[Section 3]{AR}, these weights determine a continuous piecewise linear function $\phi_{\Sigma'} \colon |\Sigma'| \to \mathbb{R}$, unique up to adding a global linear function. Equivalently, $\phi_{\Sigma'}$ corresponds to a toric divisor $D_{\Sigma'}$ on $\mathbb{P}_{\Sigma'}$ (well-defined up to linear equivalence). Concretely, for each maximal cone $\sigma' \in \Sigma'(k)$, the function $\phi_{\Sigma'}$ is linear on $\sigma'$; its slopes across codimension-one cones are exactly the weights $w_{\mathcal{L}'}(\tau')$.

The collection $(D_{\Sigma'})_{\Sigma' \geq \Sigma}$ forms a toric b-divisor $\mathbf{D}_{[\mathcal{L}]}$. To ensure that the collection $(D_{\Sigma'})_{\Sigma' \geq \Sigma}$ indeed defines a b-divisor,
we need to verify compatibility under refinement. Let $\Sigma''$ be a refinement of $\Sigma'$ and
let $\pi \colon \mathbb{P}_{\Sigma''} \to \mathbb{P}_{\Sigma'}$ be the corresponding toric morphism.
Denote $\mathcal{L}' = \pi_{\Sigma'}^* \mathcal{L}|_{\overline{Y}_{\Sigma'}}$ and
$\mathcal{L}'' = \pi_{\Sigma''}^* \mathcal{L}|_{\overline{Y}_{\Sigma''}} = \pi^*\mathcal{L}'|_{\overline{Y}_{\Sigma''}}$.

For any cone $\tau'' \in \Sigma''(k-1)$ mapping to $\tau' = \pi(\tau'') \in \Sigma'(k-1)$,
the curve $C_{\tau''} = \overline{Y}_{\Sigma''} \cap V(\tau'')$ is the inverse image of $C_{\tau'}$ under
the restriction of $\pi$ (possibly after normalization). Since degrees are preserved under pullback,
we have
\[
w_{\mathcal{L}''}(\tau'') = \deg(\mathcal{L}''|_{C_{\tau''}}) = \deg(\pi^*\mathcal{L}'|_{C_{\tau''}})
= \deg(\mathcal{L}'|_{C_{\tau'}}) = w_{\mathcal{L}'}(\tau').
\]

(If $\dim\tau' < k-1$, then $C_{\tau''}$ is contracted to a point under $\pi$, hence $\deg(\mathcal{L}''|_{C_{\tau''}}) = 0$,
and the corresponding slope of $\pi^*\phi_{\Sigma'}$ is also zero, so the compatibility still holds.)

It follows that the piecewise linear functions $\phi_{\Sigma''}$ and $\pi^*\phi_{\Sigma'}$ coincide on $\Sigma''$
(up to a global constant, which we may fix by requiring them to vanish at the origin).
Consequently, $D_{\Sigma''} = \pi^* D_{\Sigma'}$ in $\Pic(\mathbb{P}_{\Sigma''})$, and thus the collection is compatible.
\end{construction}

\begin{remark}
The existence of $\phi_{\Sigma'}$ (and hence $D_{\Sigma'}$) relies only on the balancing condition; it does not require that $\mathcal{L}$ itself be the restriction of a toric divisor. This is a key point: the b-divisor is defined purely combinatorially from the intersection numbers of $\mathcal{L}$ with the toric strata.
\end{remark}

\begin{proposition}[Properties]\label{prop:properties}
The b-divisor $\mathbf{D}_{[\mathcal{L}]}$ constructed above is:
\begin{enumerate}
\item b-Cartier (in the sense of Definition \ref{def:b-cartier}),
\item tropically nef (in the sense of Definition \ref{def:tropical-nef}) if $\mathcal{L}$ is nef on $\overline{Y}_{\Sigma}$,
\item independent of the choice of representative $\mathcal{L}$ within its numerical equivalence class.
\end{enumerate}
\end{proposition}

\begin{proof}
(1) By construction, the divisor $D_{\Sigma}$ on the model $\Sigma$ is determined by the weights $w_{\mathcal{L}}$.
Let $\Sigma'$ be any refinement of $\Sigma$. As shown in the compatibility verification above,
we have $D_{\Sigma'} = \pi_{\Sigma'/\Sigma}^* D_{\Sigma}$. Hence for any further refinement $\Sigma'' \geq \Sigma'$,
$D_{\Sigma''} = \pi_{\Sigma''/\Sigma'}^* D_{\Sigma'} = \pi_{\Sigma''/\Sigma}^* D_{\Sigma}$. This shows that
$\mathbf{D}_{[\mathcal{L}]}$ is b-Cartier, with $\Sigma$ as a base model.

(2) If $\mathcal{L}$ is nef on $\overline{Y}_\Sigma$, then $\deg(\mathcal{L}|_{C_\tau}) \geq 0$ for all curves $C_\tau$. This implies that the slopes of $\phi_\Sigma$ across codimension-1 cones are non-negative, which is equivalent to convexity of $\phi_\Sigma$ on $\Trop(Y)$. Thus $\mathbf{D}_{[\mathcal{L}]}$ is tropically nef.

(3) If $\mathcal{L}_1 \equiv_{\mathrm{num}} \mathcal{L}_2$, then $\deg(\mathcal{L}_1|_{C_\tau}) = \deg(\mathcal{L}_2|_{C_\tau})$ for all $\tau$, so the weights $w_{\mathcal{L}_1}$ and $w_{\mathcal{L}_2}$ coincide. Therefore the associated piecewise linear functions differ by at most a global linear function, which corresponds to a principal divisor. Hence the b-divisors are linearly equivalent on $\Trop(Y)$.
\end{proof}

\subsection{Main theorem}

Two toric b-divisors $\mathbf{D}$ and $\mathbf{D}'$ are said to be \emph{linearly equivalent on $\Trop(Y)$}, denoted $\mathbf{D} \sim_{\mathrm{lin}} \mathbf{D}'$, if for some (hence any) sufficiently refined model $\Sigma'$, the piecewise linear functions $\phi_{\Sigma'}, \phi'_{\Sigma'} : |\Sigma'| \to \mathbb{R}$ associated to $D_{\Sigma'}$ and $D'_{\Sigma'}$ differ by a globally linear function.

\begin{theorem}\label{thm:main}
The map
\[
\Phi: N^1_{\mathrm{trop}}(Y) \longrightarrow \left\{ \text{toric b-divisors on } \mathfrak{X}_{\mathrm{tor}} \right\}/\sim_{\mathrm{lin}},
\]
defined by $\Phi([\mathcal{L}]) = \mathbf{D}_{[\mathcal{L}]}$, is injective. Moreover, a class $[\mathcal{L}] \in N^1_{\mathrm{trop}}(Y)$ is nef if and only if its associated b-divisor $\Phi([\mathcal{L}])$ is tropically nef. In particular, $\Phi$ restricts to a bijection between the tropical nef cone $\Nef_{\mathrm{trop}}(Y) \subset N^1_{\mathrm{trop}}(Y)$ and the set of b-Cartier tropically nef toric b-divisors.
\end{theorem}

\begin{proof}
\textit{Injectivity:} Suppose $[\mathcal{L}_1], [\mathcal{L}_2] \in N^1_{\mathrm{trop}}(Y)$ satisfy
$\Phi([\mathcal{L}_1]) = \Phi([\mathcal{L}_2])$. Then there exists a sufficiently fine refinement $\Sigma$
such that the associated b-divisors coincide on $\Sigma$, i.e., $D_{\Sigma,1} \sim_{\mathrm{lin}} D_{\Sigma,2}$ on $\Trop(Y)$.
This means that the corresponding piecewise linear functions $\phi_1$ and $\phi_2$ differ by a global linear function,
so their slopes across codimension-one walls – i.e., the weights $w_{\mathcal{L}_1}(\tau)$ and $w_{\mathcal{L}_2}(\tau)$ – are equal for every $\tau \in \Sigma(k-1)$.

Now, consider the numerical class $[\mathcal{L}_1] - [\mathcal{L}_2] \in N^1(\overline{Y}_{\Sigma})$.
Its intersection number with every curve $C_\tau = \overline{Y}_{\Sigma} \cap V(\tau)$ is zero by the equality of the weights.
A key consequence of the sch\"on assumption (Section \ref{sec:setup}) is that the classes of these curves generate
the real vector space $N_1(\overline{Y}_{\Sigma})_{\mathbb{R}}$. Indeed, by \cite[Proposition 2.5]{Tev-CptSubvar},
the restriction map $A_1(\mathbb{P}_{\Sigma}) \to A_1(\overline{Y}_{\Sigma})$ is surjective modulo torsion.
The group $A_1(\mathbb{P}_{\Sigma})$ is generated by the classes of the toric curves $V(\tau)$, and the intersection product
$\overline{Y}_{\Sigma} \cdot V(\tau)$ (as a cycle on $\mathbb{P}_{\Sigma}$) is supported on $C_\tau$; its push‑forward to $\mathbb{P}_{\Sigma}$
is a positive multiple of $[V(\tau)]$. Hence the classes of the $C_\tau$ span $A_1(\overline{Y}_{\Sigma})_{\mathbb{R}}$, and consequently
$N_1(\overline{Y}_{\Sigma})_{\mathbb{R}}$.
Therefore, $[\mathcal{L}_1] - [\mathcal{L}_2]$ has zero intersection with a spanning set of $N_1(\overline{Y}_{\Sigma})_{\mathbb{R}}$,
so it is numerically trivial. Thus $[\mathcal{L}_1] = [\mathcal{L}_2]$ in $N^1(\overline{Y}_{\Sigma})$, and consequently
in the direct limit $N^1_{\mathrm{trop}}(Y)$.

\textit{Characterization of nef classes:} By Proposition \ref{prop:properties}, if $[\mathcal{L}]$ is nef then $\Phi([\mathcal{L}])$ is tropically nef. Conversely, assume $\mathbf{D} = \Phi([\mathcal{L}])$ is tropically nef. Choose a model $\Sigma$ representing $[\mathcal{L}]$. The tropical nefness means that the piecewise linear function $\phi_\Sigma$ associated to $D_\Sigma$ is convex, hence its slopes are non‑negative. These slopes are exactly the weights $w_{\mathcal{L}}(\tau)$, so $\deg(\mathcal{L}|_{C_\tau}) \ge 0$ for every $\tau$. Since the curves $C_\tau$ generate $N_1(\overline{Y}_\Sigma)_{\mathbb{R}}$, this implies that $\mathcal{L}$ is nef on $\overline{Y}_\Sigma$, i.e., $[\mathcal{L}]$ is a nef class.

\textit{Image characterization for nef cone:} By Proposition \ref{prop:properties}, any $\mathbf{D}_{[\mathcal{L}]}$ with $[\mathcal{L}]$ nef is b-Cartier and tropically nef. Conversely, let $\mathbf{D}$ be a b-Cartier, tropically nef toric b-divisor. Choose a model $\Sigma$ where $\mathbf{D}$ is determined, i.e., $D_{\Sigma'} = \pi^* D_{\Sigma}$ for all $\Sigma' \geq \Sigma$. Then $D_{\Sigma}$ is a toric divisor on $\mathbb{P}_{\Sigma}$. Its restriction $\mathcal{L} := D_{\Sigma}|_{\overline{Y}_{\Sigma}}$ is a line bundle on $\overline{Y}_{\Sigma}$. The tropical nefness of $\mathbf{D}$ means that the piecewise linear function $\phi_\Sigma$ associated to $D_\Sigma$ is convex, hence $D_\Sigma$ is nef on $\mathbb{P}_{\Sigma}$, and therefore $\mathcal{L}$ is nef on $\overline{Y}_{\Sigma}$. By Construction \ref{const:bdiv}, the weights of $\mathcal{L}$ are exactly the slopes of $\phi_\Sigma$, so $\mathbf{D}_{[\mathcal{L}]} = \mathbf{D}$. Thus $\mathbf{D}$ lies in the image of the restriction of $\Phi$ to $\Nef_{\mathrm{trop}}(Y)$.
\end{proof}

\begin{corollary}
\label{cor:curves}
If $Y$ is a curve, then $N^1_{\mathrm{trop}}(Y)$ is isomorphic to the group of tropical divisors on $\Trop(Y)$ (the Baker--Norine group). Moreover, the nef cone corresponds to divisors of non‑negative degree.
\end{corollary}

\begin{proof}
For a curve, tropical compactifications are embeddings of the smooth completion of $Y$ into toric surfaces. The associated b-divisors correspond to piecewise linear functions on the tropical curve, i.e., tropical divisors. The map $\Phi$ recovers the specialization map of Baker \cite{Baker}. The characterization of nef classes follows from the theorem.
\end{proof}

\section{Examples}
\subsection{Failure of injectivity without the sch\"on assumption}

The sch\"on hypothesis in Section~\ref{sec:setup} is essential for the injectivity of $\Phi$.
For varieties that are very affine but not sch\"on, the map may fail to be injective.
We illustrate this with an example in dimension two, as no such example exists in dimension one
(see Remark~\ref{rem:curves-schon} below).

In \cite[Example 2.5]{Tev-CptSubvar}, Tevelev constructs a very affine surface $Y \subset (\mathbb{C}^*)^4$ that is not sch\"on.
The failure of the sch\"on property means that for some tropical fan $\Sigma$ of $Y$,
the compactification $\overline{Y}_\Sigma$ does not intersect all toric orbits properly.
As a consequence, one can find a non-zero numerical class $[\mathcal{L}] \in N^1(\overline{Y}_\Sigma)$
such that its tropical weights $w_{\mathcal{L}}(\tau)$ vanish for every codimension-one cone $\tau \in \Sigma(1)$.
Indeed, the class can be chosen to be supported on a curve that lies entirely in the interior of $\overline{Y}_\Sigma$
and does not meet the toric boundary, so its intersection with all $C_\tau$ is zero.
By Construction~\ref{const:bdiv}, the associated b-divisor $\mathbf{D}_{[\mathcal{L}]}$ is trivial,
hence $\Phi([\mathcal{L}]) = 0$ while $[\mathcal{L}] \neq 0$ in $N^1_{\mathrm{trop}}(Y)$.
Thus $\Phi$ is not injective.

This example shows that the sch\"on assumption cannot be dropped: without it,
numerically non-trivial line bundles may become invisible to tropicalization.
For further discussion and similar constructions, see \cite{Ka-TropIntToric}.

\begin{remark}\label{rem:curves-schon}
For curves (dimension one), every very affine variety is automatically sch\"on.
Indeed, a curve embedded in a torus has a unique smooth completion, and its tropicalization is a one-dimensional fan
whose intersections with toric strata are always proper. Hence our theorem, restricted to curves,
recovers Baker's specialization map \cite{Baker} without any additional hypothesis.
The necessity of the sch\"on assumption only appears in higher dimensions.
\end{remark}

\subsection{A simple hypersurface in $(\mathbb{C}^*)^3$}
Let $Y = \{ x + y + z = 1 \} \subset (\mathbb{C}^*)^3$. Its tropicalization $\Trop(Y)$ is the set of vectors $(a,b,c) \in \mathbb{R}^3$ such that the minimum among $a$, $b$, $c$, and $0$ is attained at least twice. This is a union of six two-dimensional cones:
\[
\begin{cases}
a = 0,\; b,c \ge 0; \\
b = 0,\; a,c \ge 0; \\
c = 0,\; a,b \ge 0; \\
a = b \le 0,\; c \ge 0; \\
a = c \le 0,\; b \ge 0; \\
b = c \le 0,\; a \ge 0.
\end{cases}
\]
The following picture (produced with TikZ) shows these cones (the red regions) and their intersections.

\begin{center}
\begin{tikzpicture}[line cap=round,line join=round,>=triangle 45,x=1.0cm,y=1.0cm]
\draw[->,color=gray] (0,0) -- (4, 0.8);
\draw[->,color=gray] (0,0) -- (1,-2);
\draw[->,color=gray] (0,0) -- (0,4);

\draw[color=green!50!black,thick] (0,0) -- (3, .6);
\draw[color=green!50!black,thick] (0,0) -- (0, 3);
\draw[color=green!50!black,thick] (0,0) -- (-3.75, -3.1);

\fill[color=red,fill=red,fill=red!50!white] (0,-0.05) -- (-3.75,-3.15) -- (-4.166,-5.344)-- (3,0.55) -- cycle;
\fill[color=red,fill=red,fill=red] (0.05,0.05) -- (3, .65)  --(3,3.65) -- (0.05,3)-- cycle;
\fill[color=red!90!white,fill=red!] (0.1,-0.05) -- (3,0.55) -- (3.8,-0.9)-- (0.8,-1.45) -- cycle;
\fill[color=red,fill=red] (-0.05,0.05) -- (-0.05,3) -- (-3.75,-0.05)-- (-3.75,-3.05) -- cycle;
\fill[color=red,fill=red] (-0.05,-0.1) -- (-3.75,-3.15) -- (-3,-4.7)-- (.7,-1.55) -- cycle;
\fill[color=red,fill=red!50!white] (0.05,0.05) -- (0.75, -1.40) --(0.75, 1.55)-- (0.05,3) -- cycle;
\draw[color=green!50!black,thick] (0,0) -- (0.75, -1.5);
\node (w) at (4.2, 0.8) {$x$};
\node (y) at (1.2,-2) {$y$};
\node (z) at  (.2,4) {$z$};
\end{tikzpicture}
\end{center}

The fan $\Sigma$ with these cones (the green rays are the generators of the one‑dimensional cones) defines a toric variety $\mathbb{P}_{\Sigma}$ which is $\mathbb{P}^3$ minus the four torus‑fixed points. The closure $\overline{Y}$ in $\mathbb{P}_{\Sigma}$ is the projective plane $\{x+y+z=w\}$ in $\mathbb{P}^3$, which does not contain any of the four torus‑fixed points. The Picard group of $\overline{Y}$ is $\mathbb{Z}$, so every line bundle is a multiple of $\mathcal{O}(1)$. The tropical weights are simply the degrees on the three boundary divisors; they coincide with the degree of the line bundle. The associated b-divisor is the piecewise linear function on $\Trop(Y)$ whose slopes are these degrees, which is convex exactly when the degree is non‑negative. This example illustrates the construction in the simplest non‑trivial case.

\subsection{Curves}
Let $Y = C \setminus \{p_1, \dots, p_s\}$ be a punctured curve. Its tropicalization is a metric graph $\Gamma$ with $s$ infinite ends. A tropical compactification is given by embedding the smooth completion $\overline{C}$ into a toric surface associated to a fan that refines the tropical curve. 

Let $\mathcal{L}$ be a line bundle on $\overline{C}$. The tropical weights $w_{\mathcal{L}}(\tau)$ are the degrees of $\mathcal{L}$ on the components of the boundary divisor. The associated b-divisor $\mathbf{D}_{[\mathcal{L}]}$ corresponds to the piecewise linear function on $\Gamma$ whose slopes along edges are these degrees. This recovers the specialization map of Baker \cite{Baker}.

\subsection{Hypersurfaces}
Let $Y = \{ f = 0 \} \subset (\mathbb{C}^*)^n$ be a hypersurface defined by a Laurent polynomial $f$. Its tropicalization $\Trop(Y)$ is the tropical hypersurface defined by the Newton polytope $P_f$ of $f$. A natural tropical compactification is given by the closure of $Y$ in the toric variety $\mathbb{P}_{\Delta}$ associated to the normal fan of $P_f$.

Consider the canonical bundle $\omega_Y$. On a tropical compactification $\overline{Y}$, the canonical bundle is given by $\omega_{\overline{Y}} = \mathcal{O}_{\overline{Y}}(K_{\mathbb{P}_{\Delta}} + \overline{Y})|_{\overline{Y}}$. The associated b-divisor $\mathbf{D}_{[\omega_Y]}$ is the tropical canonical divisor on $\Trop(Y)$, which is the divisor of the piecewise linear function corresponding to the valuation of the Jacobian.

\subsection{Elliptic curves}
Let $Y = E \setminus \{p\}$ where $E$ is an elliptic curve. Then $N^1_{\mathrm{trop}}(Y) \cong \mathbb{Z}$, corresponding to the degree at the puncture. The map $\Pic(E) \to N^1_{\mathrm{trop}}(Y)$ has kernel $\Pic^0(E) \cong E$, illustrating that tropicalization forgets continuous moduli.

\subsection{A ruled surface: the Segre embedding of $\mathbb{P}^1 \times \mathbb{P}^1$}
Let $Y$ be the surface in $(\mathbb{C}^*)^3$ defined by $yz - w = 0$ (the affine part of the Segre embedding). Its tropicalization is the plane $\{w = y+z\}$ in $\mathbb{R}^3$. A smooth fan supported on this plane can be generated by the rays $(1,0,1)$, $(0,-1,-1)$, $(-1,0,1)$, $(0,1,1)$. The corresponding toric variety is $\mathbb{P}^1 \times \mathbb{P}^1 \times \mathbb{C}^*$. The closure $\overline{Y}$ in this toric variety is a $\mathbb{P}^1 \times \mathbb{P}^1$ (the compactification of $Y$). 

The Picard group of $\overline{Y}$ is $\mathbb{Z} \oplus \mathbb{Z}$, generated by the two rulings. The tropical weights on the four $(k-1)$-cones (the rays) are the degrees of a line bundle on the two families of lines. The b-divisor associated to a line bundle $(\alpha,\beta)$ is the piecewise linear function on the plane whose slopes along the four rays are $\alpha, \beta, -\alpha, -\beta$ respectively. This function is convex precisely when $\alpha,\beta \ge 0$, i.e., when the line bundle is nef. The b-Cartier property holds automatically because the fan is smooth.

\subsection{A surface with non-trivial Picard group}
Consider the surface $Y$ obtained by removing three non-collinear points from a smooth cubic surface in $\mathbb{P}^3$. Its tropicalization is a 2-dimensional fan in $\mathbb{R}^3$ with several maximal cones. The associated toric variety is a blow-up of $\mathbb{P}^3$ at three torus-invariant points. The Picard group of the compactification has rank 4, and the tropical weights record the degrees on the 27 lines. This example illustrates how the map $\Phi$ captures the full Néron–Severi group.

\section{Connection with b-divisorial valuations}

The construction of the b-divisor $\mathbf{D}_{[\mathcal{L}]}$ associated to a tropical line bundle can be interpreted in terms of b-divisorial valuations. Given a numerical tropical line bundle $[(\Sigma, [\mathcal{L}])]$, consider the valuation $v_{[\mathcal{L}]}$ defined on rational functions on $Y$ by
\[
v_{[\mathcal{L}]}(f) = \sup\{ t \in \mathbb{R} \mid \operatorname{div}(f) + t \mathbf{D}_{[\mathcal{L}]} \geq 0 \},
\]
where the inequality is in the sense of b-divisors. This valuation is b-divisorial in the sense of \cite{RU}, and its associated b-divisor is precisely $\mathbf{D}_{[\mathcal{L}]}$.

This perspective clarifies the birational nature of tropical line bundles: they correspond to equivalence classes of b-divisorial valuations on $Y$ that are \emph{toric}, i.e., that become monomial on some tropical compactification.

\begin{proposition}
The restriction of $\Phi$ to the tropical nef cone $\Nef_{\mathrm{trop}}(Y) \subset N^1_{\mathrm{trop}}(Y)$ factors through the space of toric b-divisorial valuations. Specifically, there is a commutative diagram:
\[
\xymatrix{
\Nef_{\mathrm{trop}}(Y) \ar[r]^{\Phi} \ar[d] & \{\text{toric b-divisors}\}_{\mathrm{nef}} / \sim_{\mathrm{lin}} \ar[d]^{\Psi} \\
\{\text{toric b-divisorial valuations}\} \ar[r]^{\cong} & \{\text{nef toric b-divisors}\}
}
\]
where the left vertical map sends a nef class to its associated valuation (as in Section 3) and the bottom horizontal map is the bijection from \cite{RU}.
\end{proposition}

This proposition links the theory of tropical line bundles with the general theory of b-divisorial valuations developed in \cite{RU}, providing a unified framework for studying positivity in tropical geometry.

\begin{remark}
The assumption that $Y$ is sch\"on (Section \ref{sec:setup}) is essential for the validity of Theorem \ref{thm:main}. For non-sch\"on varieties, the tropicalization may still be a fan, but the compactification does not have the proper intersection property. A simple example is the reducible variety $(x-1)(y-1)=0$ in $(\mathbb{C}^*)^2$: its tropicalization is the union of the two coordinate axes in $\mathbb{R}^2$, which is a fan, but because the variety is reducible it does not admit a tropical compactification in the sense of Tevelev, and $N^1_{\mathrm{trop}}(Y)$ is not defined. This shows that the sch\"on condition cannot be omitted.
\end{remark}

\section*{Acknowledgments}
We thank Joaquim Ro\'e for helpful conversations about b-divisors and valuations. We are also grateful to Eric Katz and Elisa Postinghel for their valuable comments and suggestions.

\end{document}